\documentclass[12pt]{article}
\title{Brownian beads}
\author{B\'alint Vir\'ag\footnote{Department of Mathematics, MIT, Cambridge, MA
02139, USA. {\sf balint@math.mit.edu}. Research partially
supported by NSF grant \#DMS-0206781.}}
\date{May 27, 2003}

\usepackage{amsfonts}
\usepackage{graphicx}
\usepackage{amsmath}
\usepackage{amsthm}

\newtheorem{theorem}{Theorem}
\newtheorem{definition}[theorem]{Definition}
\newtheorem{question}[theorem]{Question}
\newtheorem{conjecture}[theorem]{Conjecture}
\newtheorem{fact}[theorem]{Fact}
\newtheorem{remark}[theorem]{Remark}
\newtheorem{proposition}[theorem]{Proposition}
\newtheorem{lemma}[theorem]{Lemma}
\newtheorem{corollary}[theorem]{Corollary}
\newenvironment{theorem*}[1]{\par \trivlist
 \itemindent 0pt \item[\hskip\labelsep\bf Theorem #1]
 \it\ignorespaces}{\endtrivlist}

\renewenvironment{proof}{\par \trivlist
 \itemindent\parindent \item[\hskip\labelsep\sc Proof.]
 \ignorespaces}{\qed\endtrivlist}
\newenvironment{proofof}[1]{\par \trivlist
 \itemindent\parindent \item[\hskip\labelsep\sc Proof of #1.]
 \ignorespaces}{\qed\endtrivlist}

\newcommand\mnote[1]{}
\newcommand\lb[1]{\label{#1}\mnote{#1}}
\newcommand\bel[1]{{\mnote{#1}}\begin{equation}\label{#1}}
\newcommand\ee{\end{equation}}

\newcommand\re[1]{(\ref{#1})}
\newcommand{\Strip}{\mathbb S}
\newcommand{\eps}{\varepsilon}

\newcommand{\Z}{\mathbb Z}

\newcommand{\RR}{\mathbb R}
\newcommand{\nis}{\hbox{$\cap$\kern-.6em\lower0em\hbox{$/$}}\,}

\newcommand{\CC}{\mathbb C}

\newcommand{\HH}{{\mathbb H}}

\newcommand{\FF}{\mathcal F}
\newcommand{\GG}{\mathcal G}

\newcommand{\ev}{{\mathbf E}}
\newcommand{\pr}{{\mathbf P}}
\newcommand{\one}{{\mathbf 1}}
\newcommand{\gi}{\,|\,}

\newcommand{\dist}{\mbox{dist}}
\newcommand{\Ghat}{{\hat G}}

\newcommand{\BE}{\mathsf{BE}}
\newcommand{\bm}{\mathsf{BM}}
\newcommand{\BB}{\mathsf{BB}}
\newcommand{\nbb}{\mathsf{NBB}}

\newcommand{\cd}{\Rightarrow}
\newcommand{\omu}{{\overline \mu}}

\newcommand{\diam}{\operatorname{diam}}
\newcommand{\capz}{\mbox{\rm cap}_0}
\newcommand{\capo}{\mbox{\rm cap}_1}

\usepackage{natbib}
\begin{document}
\maketitle
\begin{abstract} \noindent
We show that the past and future of half-plane Brownian motion
at certain cutpoints are independent of each other after a
conformal transformation. Like in It\^o's excursion theory, the
pieces between cutpoints form a Poisson process with respect to
a local time. The size of the path as a function of this local
time is a stable subordinator whose index is given by the
exponent of the probability that a stretch of the path has no
cutpoint. The index is computed and equals $1/2$.
\end{abstract}

\section{Introduction}

The mathematical theory of conformally invariant planar stochastic
models has seen great progress in the recent years. The goal of
this paper is to consider Brownian motion, the first example of a
conformally invariant process, and further explore its conformal
structure.

Let $B$ denote Brownian motion started at zero and conditioned to
stay in the upper half plane; we call this distribution $\BE$, or
half-plane excursion. This is a transient process, and he path has
cut-points \citep{burdzy89}, that is points which, if removed,
make the image of $B$ disconnected. We call the segments of the
paths between consecutive cutpoints {\bf Brownian beads}.

For a given cutpoint, the complement of the past path has one
infinite connected component. This can be mapped back to the half
plane via a conformal homeomorphism. It is convenient to normalize
this map so that it has derivative 1 at $\infty$ and takes the
cutpoint to 0; we call this map the {\bf conformal shift}.
\begin{figure}
\label{f2} \centering
\includegraphics[height=3in]{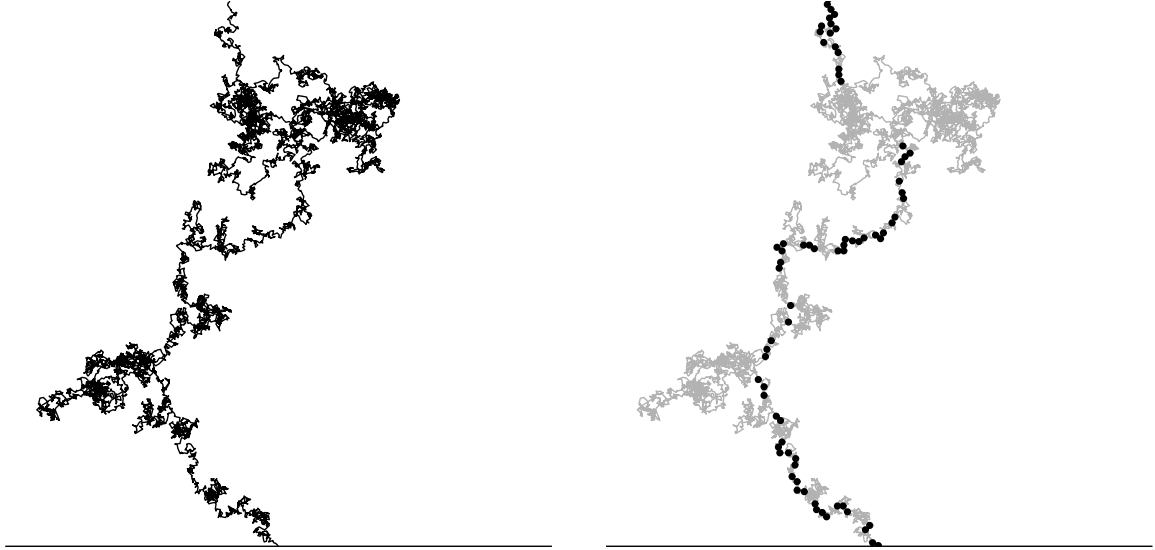}
\caption{Brownian excursion in $\HH$ and its cutpoints}
\end{figure}
Our first goal is to show that Brownian beads are independent
of each other after a conformal transformation. More precisely,
for a cutpoint $g$ let $f$ denote the corresponding conformal
shift and let $\beta$ the bead starting at $g$ (if any).
\begin{theorem}\lb{ppp}
There exists a local time $\lambda$ supported on cutpoints so that
 $(f_{\lambda}(\beta_{\lambda}), \lambda \ge 0)$ is a Poisson
point process.
\end{theorem}
Here $f_\lambda$ applied to a path maps the image of the path
according to the conformal map, and changes time-parameterization
according to local Brownian scaling (Section \ref{s.conformal}).
The intensity measure of the Poisson point process is a
$\sigma$-finite measure on paths conditioned to start without
cutpoints; we call this measure the bead process. In Section
\ref{s.properties} we prove a Markov property for the bead process
and show that it has finite lifetime.

Theorem \ref{ppp} may be thought of as a two-dimensional analogy
of It\^o's theorem about Brownian excursions. Recall that
excursions of 1-dimensional Brownian motion are the segments of
the path between two consecutive visits to 0. For a zero $g'$ let
$\beta'$ denote the excursion starting at $g'$, and let $f'$
denote the map that shifts the future of the path by $-g'$.
\begin{theorem}[It\^o]
There exists a local time $\lambda'$ on zeros so that
$(f'_{\lambda'}(\beta'_{\lambda'}), \lambda' \ge 0)$ is a Poisson
point process.
\end{theorem}
The local time $\lambda'$ at zero has the nice property that the
$(g'(\lambda'),\lambda'\ge 0)$ process, or in other words the
total duration of excursions up to the zero $g'$ is a stable
subordinator of index $1/2$. It turns out that the analogous
statement holds in our setting. Let $a(t)$ denote the half-plane
capacity of the path up to time $t$; half-plane capacity is a way
to measure the size of subsets of $\HH$ (see Section
\ref{s.conformal}). Then
\begin{theorem}
The process $(a\circ g(\lambda), \lambda \ge 0)$ is a stable
subordinator of index $\alpha$.
\end{theorem}
Another property of beads which is shared by  It\^o excursions is
that the process $(\beta_\lambda,\,\lambda\ge 0)$ determines
$(B(t),\,t\ge0)$ (see Remark \ref{determine}).

In Section \ref{s.exponent} we show that the index $\alpha$ can be
expressed as an exponent.

\begin{theorem}\lb{exponent} We have $\pr(\mbox{B has no cuttime in }(1,t))
= t^{-\alpha + o(1)}.$
\end{theorem}

Finally, we compute $\alpha$.
 \newcommand{\indexhalfthm}
 {The Brownian bead index $\alpha$ equals $1/2$.}
 \begin{theorem}\lb{indexhalf}
 \indexhalfthm
 \end{theorem}
Although the $1/2$ here is equal to the index for It\^o's
excursions, it is a coincidence. Neither is it a consequence of
Brownian scaling. Unlike the proof of most exponents, this theorem
does not rely on SLE processes (nor exponents that use these
processes). The proof uses beads, not the exponent representation;
I am grateful to W. Werner with whom I discussed ideas of this
kind. A related exponent governing the Hausdorff dimension of
cuttimes for $\mbox{SLE}_6$ was computed by \cite{beffara}, and it
seems possible to derive the value of $\alpha$ directly from his
results. \cite{math.PR/0302115} also computed an exponent related
to $\alpha$ using SLE$(\kappa,\rho)$ processes (in a paper with
beautiful ideas and computations but few rigorous proofs). The
advantage of the proof here is that it is more elementary and
conceptual; it relies only on the special intersection exponent
$\xi(2,1)$, whose value was determined by \cite{lawler95}. His
proof covers even the 3-dimensional case.

In the proof, we will use the following fact about Brownian
excursion. Let $A$ be a compact subset of $\overline \HH\setminus
\{0\}$ so that the sets $\RR \cup A$ and $\HH\setminus A$ are
connected. Let $f:\HH\setminus A\to\HH$ be a conformal
homeomorphism fixing $\infty$ with $f'(\infty)=1$.
\begin{proposition}\lb{mainhit}
$\quad\pr(B\mbox{ avoids } A) =f'(0)$.
\end{proposition}

Further interesting properties of Brownian motion have recently
been found using the fact that its outer boundary has the same
distribution as (the outer boundary of) some SLE processes, see
\cite*{lwsrest} and the references therein. Some of methods used
in this paper appear independently in the recent literature; in
particular, see \cite{dubedat} and \cite{loopsoup}.

Many of the arguments of this paper carry over to ordinary
Brownian motion, more precisely, Brownian excursion from a
boundary point to an interior point of a domain. The analysis
there is a bit more difficult since there is no scale invariance;
we do not follow this avenue here.

In Section \ref{properties} we discuss some simple properties of
Brownian excursions in planar domains that we will need later. We
also show Proposition \ref{mainhit} there. In Section \ref{s.ctf}
we introduce the cuttime filtration and prove a Markov property
with respect to this filtration. Section \ref{s.conformal} reviews
the facts we need from conformal geometry, and introduces some
basic semigroups of paths. In Section \ref{s.independent} we prove
the Poisson point process decomposition of Theorem \ref{ppp}, and
we define the bead process, a Brownian excursion in the half plane
conditioned to start without cutpoints. In the next section we
show that this special process does not need infinite time to
start off, a fact that is not clear from the definition. In the
last two sections, \ref{s.exponent} and \ref{s.index}, we prove
Theorems \ref{exponent} and \ref{indexhalf}; we also give some
open questions and conjectures.

\section{Properties of Brownian excursion}\lb{properties}

Let $D$ be a {\bf regular domain}, i.e. a simply connected open
subset of the plane whose boundary is locally connected. Let
$\partial D$ denote the Caratheodory boundary, i.e. the set of
prime ends of $D$.

For $\{a\},Z\subset D\cup \partial D$, and $a\not=\infty$, let
$\BE(a,Z,D)$ denote Brownian motion started at $a$ and conditioned
to hit $Z$ no later than $\partial D$. If this event has positive
probability, then this definition is precise; otherwise it can be
made precise by considering $h$-processes or by taking a limit. A
special case is $\BE(0,\infty,\HH)$, the half-plane excursion,
which we will often abbreviate $\BE$. The coordinates of $\BE$ are
Brownian motion and an independent 3-dimensional Bessel process.

The two most important properties of $B\sim\BE(a,Z,D)$ are
restriction and conformal invariance. Restriction says that if
$D'\subset D$, then on the event that $B$ stays in $D'$,  it
has distribution $\BE(a,Z',D')$, where $Z'=Z\cap(D'\cup
\partial D')$. Conformal invariance says that if
$f:(a,Z,D)\to (a',Z',D')$ is a conformal homeomorphism between two
regular domains, then after a time change (\ref{timechange}) the
image of $B$ under $f$ has distribution $\BE(a',Z',D')$. The
restriction property is straightforward; conformal invariance has
been first proved by \cite{levy40}; the later development of
stochastic calculus makes it a simple exercise.

It is important for the previous paragraph that by Theorem 9.8
in the book of \cite{pommerenke}, conformal homeomorphisms from
the half plane $\HH$ or unit disk to regular domains extend
continuously to the boundary.

\begin{remark}\lb{curve} \rm
One way to get such domains is the following. Let $E$ denote the
image of a curve $\gamma:[0,1]\rightarrow \CC$, and let $D$ be a
connected component of $\CC\setminus E$. Since local connectivity
is preserved under continuous maps, $E$ is locally connected. The
proof of Theorem 9.8 in \cite{pommerenke}, part (ii) $\rightarrow$
(iii) carries through without changes to show that in fact the
boundary of $D$ is also locally connected.
\end{remark}

The strong Markov property says that the future of $B$ after a
stopping time $\tau$ given $\FF_\tau$ has distribution
$\BE(B(\tau),Z,D)$.

The following theorem is due to \cite{cm83} (see
\cite{MR86d:60088} for a simpler proof).

\begin{theorem}\label{cm}
There exists a universal constant $c$ so that the lifetime $\tau$
of a path distributed $\BE(a,Z,D)$ has $\ev\tau < c\;$area$(D)$.
\end{theorem}

We say that a sequence $a_n\in D$ {\bf converges to $a\in \partial
D$ along a path} if there is a continuous curve $\gamma:(0,1)\to
D$ so that $a_n=\gamma(t_n)$ for some sequence $t_n$, and $\gamma$
is contained in the equivalence class that defines the prime end
$a$.

In order to define a distance, we may extend paths defined on a
finite interval $[0,t]$ to be constant on $(-\infty,0]$ and
$[t,\infty)$. Let
\begin{eqnarray*}
\dist_0(\pi_1,\pi_2)&=&\sup_{t\in \RR}\|\pi_1(t)-\pi_2(t)\|,\\
\dist(\pi_1,\pi_2)&=&\inf_{\eps \in
\RR}\left(|\eps|+\dist_0(\theta_\eps\pi_1,\pi_2)\right).
\end{eqnarray*}
where $\theta_\eps$ denotes time shift. Let ``$\Rightarrow$''
denote convergence in distribution (in the case of paths with
respect to the metric ``$\dist$'').

\begin{proposition}\lb{kezdet}
If $a_n\to a\in \partial D$ along a path, then $\BE(a_n,z,D)
\Rightarrow \BE(a,z,D).$
\end{proposition}

\newcommand{\cdist}{\mbox{cdist}}
\begin{proof}
Consider the conformal homeomorphism $g:\HH \rightarrow D$ which
maps $0$ to the prime end $a$ (as determined by the path of
convergence of $a_n$). Since the domain $D$ is regular, we have $
g^{-1}(a_n)\rightarrow 0 $. Let
$$N_\eps=\{g(z)\ : \ z\in \HH, |z|<\eps\}.$$
Since the map $g$ is uniformly continuous in a neighborhood of
$0$, it follows that diam$(N_\eps)\rightarrow 0$ as $\eps \to 0$.

If $a_n\in N_\eps$, then let $p_{n,\eps}$ be the probability that
independent paths with distributions $\BE(a_n,Z,D)$ and
$\BE(a,Z,D)$ do not intersect before their respective hitting
times $\tau_{n,\eps}$, $\tau_\eps$ of $N_\eps$. By conformal
invariance, this probability can be computed in the half-plane
image. Then it is well known that there exists $c,\gamma>0$ so
that
$$
p_{n,\eps}:=\| \mu_{n,\eps} - \mu_{\eps}\|_{TV} <
c(|g^{-1}(a)-g^{-1}(a_n)|/\eps)^\gamma.
$$
Indeed, one can consider times of hitting concentric circles of
radii $2^{-k}$. Between such times the probability that two
Brownian excursions intersect is bounded away from $0$. Thus
$p_{n,\eps}\rightarrow 0$ for $\eps$ fixed and $n\rightarrow
\infty$.

Let $B$ have distribution $\BE(a,z,D)$, and define $B_n$ as having
distribution $\BE(a_n,z,D)$ and coupled to $B$ at its first
hitting time of the path of $B$. The event $C_n$ that this can be
done before either process has left $N_\eps$ has probability
$1-p_{n,\eps}$.

Let $A_n$ be the event that $|\tau_{\infty,\eps}-\tau_{n,
\eps}|<\delta$. By Theorem \ref{cm} and Markov's inequality,
$\pr(A_n^c)<c_2\; $area$(N_{\eps})/\delta$, for some universal
constant $c_2$. Now if $A_n$ and $C_n$ hold, and $\eps<\delta$,
then $\dist(B,B_n)<3\delta$. Therefore
$$
\pr(\dist(B,B_n)> 3 \delta) \le c_2\;\mbox{area}(N_{\eps})/\delta
+ p_{n,\eps} + \one(\delta\le \eps)
$$
and the right hand side can be made arbitrarily small by first
picking $\eps$ then letting  $n\rightarrow \infty$. It follows
that $B_n$ converges to $B$ in probability and so in distribution.
\end{proof}

\begin{corollary}\label{01H}
Let $Z_n\rightarrow {1}$ be sets of reals, and let $a_n\in \HH$,
$a_n \rightarrow 0$. Then $$\BE(a_n,Z_n,\HH) \Rightarrow
\BE(0,1,\HH).$$
\end{corollary}

\begin{proof}
It suffices to prove this for the case when each $Z_n=\{z_n\}$ is
a single point, since the general case is a mixture of such
processes. Consider the rescaled versions $\BE(a_n/z_n,1,\HH)$ of
$\BE(a_n,z_n,\HH)$. These converge to $\BE(0,1,\HH)$ by the lemma,
and since the scaling factor converges to 1, so do the original
processes.
\end{proof}

Let $\BE(a,z,D,K)$ denote the distribution of $\BE(a,z,D)$
conditioned to hit $K$. Say a sequence of sets $A_n\subset D$ {\bf
converges to $a$ along a path} if every sequence of points $a_n\in
A_n$ does.

\begin{lemma}\label{slimhit}
Let $D$ be a regular domain, and consider regular subdomains $D_n
= D \setminus K_n$, where the $K_n \subset D$ are relatively
closed. Let $L_n \subset D_n$ closed, and let $a_n\in D_n$.
Suppose that
$$M_n:=K_n\cup L_n \cup \{a_n\} \rightarrow \{a\}\subset
\partial D  \quad \mbox{along a path}.$$ Then
$$\BE(a_n,z,D_n,L_n) \cd
\BE(a,z,D).$$
\end{lemma}

\begin{proof}
Let $W_n$ have distribution $\BE(a_n,z,D_n,L_n)$. The slightly
difficult part is to construct relatively closed sets $S_n
\rightarrow \{a\}$ separating $z$ from $M_n$ in $D$ so that
 \bel{e1}
 \pr(W_n \mbox{\ hits\ } S_n \mbox{ before }L_n)\rightarrow 0,
 \ee
 \vskip-1em
 \bel{e2}
 \sup_{s\in S_n} \pr(\BE(s,z,D) \mbox{ hits } K_n) \rightarrow 0.
 \ee

To complete the proof from here, let $\tau_n$ be the hitting time
of $S_n$ for $W_n$. Since $S_n \rightarrow \{a\}$, we have $\tau_n
\rightarrow 0$. If we condition $W_n$ to have done its duty of
hitting $L_n$ by time $\tau_n$, then $\{W_n(\tau_n+t),\ t\ge 0\}$
has the same distribution as $\BE(W(\tau_n),z,D_n)$. This, in turn
has the same distribution as $\BE(W(\tau_n),z,D)$ conditioned not
to hit $K_n$. In both cases, we are conditioning on events whose
probabilities converge to 1. Also, $\BE(W(\tau_n),z,D)$ is a
mixture of processes starting from points of $S_n$, so by
Proposition \ref{kezdet} it converges to $\BE(a,z,D)$. It follows
that the distribution of $W_n$ converges to $\BE(a,z,D)$, as
required.

Now we proceed to find the sets $S_n$. By conformal invariance
of the probabilities involved, it suffices to do this in the
upper half plane with $a=0,\ z=\infty$. Here we use that our
domain has a curve boundary, so if in $\HH$ we have $S_n
\rightarrow 0$, then for its image under the conformal
homeomorphism $(0,\infty,\HH) \to (a,z,D)$ we have $S'_n
\rightarrow a$.

Let $c_n\rightarrow \infty$ slow enough that we still have $c_n
M_n \rightarrow {0}$; if we construct $S_n$ that are bounded for
this  rescaled problem, then scaling back will make $S_n
\rightarrow {0}$. So it suffices to construct bounded $S_n$. We
consider the uniformizing map $g_n$ which takes $D_n$ to the half
plane $\HH$. Assume that $g_n$ has hydrodynamic normalization, and
extend it to the boundary of $D_n$.  The half-plane version of the
Caratheodory Kernel Theorem (see Theorem 1.8 in \cite{pommerenke})
implies that $g_n(z)\rightarrow z$ uniformly where it is defined.
So if we $L_n'=g_n(L_n)$, $a_n'=g_n(a_n)$ then
$L'_n\cup\{a'_n\}\rightarrow \{0\}$.

Let $S$ be the unit semicircle in $\HH$, and use the shorthand
$B_n=\BE(a'_n,\infty,\HH)$. We will show that
 \bel{e5}
 \pr(B_n \mbox{ hits } S \mbox{ before } L'_n\ |\ B_n \mbox{ hits }L'_n
)\rightarrow 0,
 \ee
\vskip-1em
 \bel{e6}
 \sup_{s\in S}\pr (\BE(s,\infty,\HH) \mbox{\ hits\ } L'_n) \rightarrow 0.
 \ee
Once we have this, we take $S_n= g_n^{-1}(S)$ so by uniform
convergence, $S_n \rightarrow S$. Thus  $\{S_n\}$ is bounded, and
$(\ref{e1}, \ref{e2})$ follow from conformal invariance.

The second statement (\ref{e6}) is obvious, and the first follows
since we can easily find semicircles $R_n \rightarrow \{0\}$ so
that
 \bel{e7}
 \sup_{w\in R_n}
 \pr(\BE(w,\infty,\HH)\mbox{\ hits\ } L'_n) <
 \pr(\BE(a'_n,\infty,\HH)\mbox{\ hits\ } L'_n).
 \ee
The left hand side of (\ref{e5}) can be written as
 \bel{e3}
 \pr(\mbox{$B_n$ hits $S$ before $L'_n$ and hits $L'_n$})
 \ /\ \pr(\mbox{$B_n$ hits }L'_n).
 \ee
 Using the Markov property and the topology of the setup the
numerator can be bounded above by
\begin{equation} \label{e4}
\sup_{w \in S} \pr(\BE(w,\infty,\HH)\mbox{ hits }R_n)
  \sup_{w\in R_n}\pr(\BE(w,\infty,\HH)\mbox{ hits }L'_n)
\end{equation}
and so by (\ref{e7}), expression (\ref{e3}) is bounded above by
the first factor of (\ref{e4}), which converges to 0, proving
(\ref{e5}).
\end{proof}

\begin{lemma}\label{fathit}
Let $D$ be a regular domain, and let $K\subset D$ be a relatively
compact connected set containing at least two points. Let $z\in
\partial D$. Suppose that $K_n\subset D$ are relatively compact
and converge to $K$. Suppose that $a_n \rightarrow a$. Then
$$\BE(a,z,D,K_n)\cd \BE(a_n,z,D,K).$$
\end{lemma}

\begin{proof}
Let $g$ be a bounded continuous test function for paths, and let
$B_n,\,B$ be distributed as $\BE(a_n,z,D)$ and $\BE(a,z,D)$
respectively.

Note that for $\BE(a,z,D)$-almost all paths $\pi$ (1) the function
$f_n = \one(\pi\ hits\ K_n)$  converges to $f = \one(\pi\ hits\
K)$ and (2) $f$ is continuous at $\pi$. Thus by Proposition
\ref{kezdet} we have
\begin{eqnarray*}\ev(f_n(B_n)g(B_n)) &\rightarrow& \ev(f(B)g(B)),
\\
\ev f_n(B_n) &\rightarrow& \ev f(B).
\end{eqnarray*}
Note that the hypothesis of the lemma ensures that $\ev f(B)$ is
positive. Taking the ratio of the previous two limit statements we
get the desired result.
\end{proof}

\subsection*{Hitting probabilities}

Let $D$ be a domain, and let $a,z$ be points on $\partial D
\setminus A$ in the neighborhood of which the boundary is a
differentiable curve. Let $A$ be a hull in $\overline D$ not
containing $a,z$. We show the following equivalent version of
Proposition \ref{mainhit}.

\begin{proposition}\lb{hit}
Let $f$ be a conformal homeomorphism that takes $D\setminus A$ to
$D$ and fixes $a,z$. Then
$$
\pr(\BE(a,z,D)\mbox{ avoids }A) = f'(a)f'(z).
$$
\end{proposition}

\begin{proof}
Because of conformal invariance, it is sufficient to prove this in
the half-plane $\HH $ with $a=0$, $z=1$.

The indicator function of the event above is continuous at almost
every path with respect to the excursion measure. Let $a_n=i/n$,
and $Z_n=(1,1+1/n)\subset \RR$. It suffices to show that as
$n\rightarrow \infty$ we have
\begin{equation}\label{ff1}
\pr(\BE(a_n,Z_n,\HH)\mbox{ avoids }A) \rightarrow f'(0)f'(1)
\end{equation}
since the measures on the left converge to the excursion measure
by Corollary \ref{01H}.
\begin{eqnarray}
\pr(\bm(a)\mbox{ hits }Z_n\mbox{  no later than }\partial
D)=h_{1/n}(i/n) = \nu n^{-2} + o (n^{-2}) \label{hit1}
\end{eqnarray}
Here $h_\eps$ is the harmonic function determined by the boundary
conditions $1$ on $(1,1+\eps)$ and $0$ elsewhere, and $\nu$ is
some fixed constant. The functions $h$ can be easily given
explicitly, but the formulas are unimportant, all we need is the
intuitive fact that as $\eps,z \rightarrow 0$ we have
$$h_\eps(z) = \nu\eps\Im z + o(\eps|z|)
$$
giving the second equality in (\ref{hit1}). Similarly, after a
conformal transformation by $f$
\begin{equation}\label{hit2}
\pr(\bm(a_n)\mbox{ hits }Z_n \mbox{ no later than }A\cup\partial
D)
\end{equation}
becomes
\begin{eqnarray}
\pr(\bm(f(a_n))\mbox{ hits }f(Z_n) \mbox{ no later than
}f(A\cup\partial \HH)) \quad  \nonumber
\\=h_{f(1+\eps)-1}(f(i\eps))  = f'(0)f'(1)\nu
n^{-2} + o(n^{-2}) \nonumber
\end{eqnarray}
but the ratio of (\ref{hit2}) and (\ref{hit1}) gives (\ref{ff1}),
proving the claim.
\end{proof}

\section{The cuttime filtration} \lb{s.ctf}

Consider the excursion $\{B_t\}$ with distribution $\BE(a,z,D)$
where $D$ is a regular domain, $a,z\in \partial D$. Let $G$ denote
the set of {\bf cuttimes}, that is, the set of times for which
images of the future and past are disjoint; \cite{burdzy89} showed
that such times exist with probability 1. It is easy to check that
the set $G$ is measurable and it is a closed set a.s.

It is clear that one cannot decide what the cuttimes are up to
time $t$ by only looking at the past of the process $\{B_t\}$. $G
\cap [0,t]$ is therefore not measurable with respect to the
standard filtration ${\mathcal F}_t$ of $\{B_t\}$. We therefore
introduce the {\bf cuttime filtration} $\mathcal G_t$ generated by
${\mathcal F}_t$ and $G \cap [0,t]$.

This enlargement of filtration may look large, but in fact it is
not in the following sense. Let $\Ghat_t$ denote the set of
cuttimes of the process $\{B_t\}$ restricted to the interval
$[0,t]$. If $\tau$ is the earliest time for which the set of
points $B[0,\tau]$ is revisited after time $t$, then we clearly
have
$$
G \cap [0,t] =  \Ghat_t \cap [0,\tau].
$$
Since $\Ghat_t$ has Hausdorff dimension less than one
\citep{lawler96}, and the distribution of $B_\tau$ is absolutely
continuous with respect to harmonic measure on $B[0,t]$, we have
by Makarov's theorem that  $\tau \notin \Ghat_t$ ${\mathcal
F}_t$-a.s. In particular, $G\cap[0,t]$ is determined by the
connected component $(g_t,g_t')$ of $[0,t] \setminus \Ghat_t$
containing $\tau$. The left endpoint of this interval, $g_t$, is
the last cuttime up to time $t$. Then given ${\mathcal F}_t$,
$g_t$ determines $\mathcal G_t$; its distribution  is concentrated
on the countable set of points in $\Ghat_t$ that are isolated from
the right.

Let $\rho B_s$ denote the connected component of $D\setminus
B[0,s]$ having $z$ as an interior or boundary point. By Remark
\ref{curve} if $D$ is a regular domain, then so is $\rho B_s$ for
every $s$. Let $\BE(a,z,D,K)$ denote the distribution of
$\BE(a,z,D)$ conditioned to hit $K$.

\begin{proposition}[Markov property for cuttime filtration]
\lb{mp} Let $\tau$ be a ${\mathcal G}_t$-stopping time. Then the
distribution of $\{B_{\tau+t}, t\ge 0\}$ given $\mathcal G_\tau$
is $\BE(B_t,z,\rho B_{g_t}, B[g_t,g_t'])$.
\end{proposition}

Note that, while for fixed time $\tau$, $B_{\tau}$ is an interior
point of $\rho B_{g_\tau}$ with probability one, there exists
stopping times for which $g_\tau=\tau$ a.s. An example is the
first cuttime after time 1. In these cases the ``hitting''
condition of the proposition holds trivially.

\begin{proofof}{Proposition \ref{mp}} \

Recall that $\GG_t$ is generated by $ \FF_t$ and $g_t$. Since the
distribution of $\{B_{t+s},s \ge 0\}$  given $\FF_t$ is
$\BE(B_t,z,D)$ by the Markov property, further conditioning on the
value of $g_t$ means that the future of the process $\BE(B_t,z,D)$
has to stay in $\rho B_{g_{\tau}}$, and hit
$B[g_{\tau},g'_{\tau}]$. This completes the proof for fixed time
$t$. Denote $L_t$ the set $B[g_{t},g'_{t}]$, and denote $\ev_t$
the distribution conditioned as above, that is $\BE(B_t,z,\rho
B_t, L_t)$.

What we  showed is equivalent to the following. If $X$ is a
bounded random variable on paths, and the time shift operator
$\theta_t$ is defined as $ (\theta_t B)(s)=B(s+t)$, then
$$
\ev (X\circ \theta_t \ |\ \GG_t) = \ev_t X.
$$
Thus, if $\tau$ is a $\GG_t$-stopping time concentrated on a
discrete set of values $\mathcal T$, then
 \begin{eqnarray*}
 \ev (X\circ\theta_\tau\ |\ \GG_\tau)
 &=& \sum_{t\in \mathcal T} \ev (X\circ \theta_\tau\ \one(\tau=t) \ |\ \GG_\tau)
 \ =\ \sum_{t\in \mathcal T} \ev (X\circ \theta_t\ \one(\tau=t) \ |\ \GG_t)
 \\&=& \sum_{t\in \mathcal T} \one(\tau=t) \ev (X\circ \theta_t\ \
 |\
 \GG_t)
 \ =\  \sum_{t\in \mathcal T} \one(\tau=t) \ev_t X
 \ =\ \ev_\tau X,
 \end{eqnarray*}
which completes the proof for this case. For an arbitrary
$\GG_t$-stopping time $\tau$ we take a discrete approximation by
letting $\tau_n$ be the first element of $ 2^{-n}\Z$ which is at
least $\tau$. Since the decision to stop at time $\tau_n$ is made
by time $\tau$, we have
$$
\ev(X\circ \theta_{\tau_n} \ |\ \GG_{\tau_n}) = \ev(X\circ
\theta_{\tau_n} \ |\ \GG_\tau).
$$
Since $\theta_t$ is a.s. continuous as $t \downarrow 0$, and $X$
is continuous, we have $X \circ \theta_{\tau_n} \rightarrow
X\circ\theta_{\tau}$ almost surely, and by the bounded convergence
theorem
$$
\ev(X\circ\theta_{\tau_n} \ | \ \GG_\tau) \rightarrow \ev(X\circ
\theta_\tau \ |\ \GG_\tau).
$$
Therefore it suffices to prove that $\FF_\infty$-a.s.
$\ev_{\tau_n} X \rightarrow \ev_\tau X$, or
 \bel{cd}
 \BE(B_{\tau_n}, z, \rho B_{g_{\tau_n}},L_{\tau_n}) \cd
 \BE(B_{\tau},   z, \rho B_{g_\tau},L_{\tau}).
 \ee
First we fix the path $B$ and examine what happens to $g_{\tau_n}$
and $g'_{\tau_n}$ in the limit. The definition $g_\tau= \sup (G
\cap [0,\tau])$ implies that $g_\tau$ is increasing in $\tau$ and
continuous from the right. Assume that $\tau\not=g_\tau$. Then
$\tau$ lies in a connected component of $G$ whose left endpoint is
$g_\tau$, so for all large $n>N_1$ we must have $g_{\tau_n}=
g_\tau$.

Now assume also that $g'_\tau \not=\tau$. Then $B[g_\tau,
g_\tau']$ is disjoint from $B(g_\tau', \tau]$. Let $t$ be the
first time after $\tau$ that $B_t$ revisits the set
$B[g_\tau,g_\tau']$. Then for $\tau_n<t$ and $n>N_1$ we have
$g'_{\tau_n}=g'_\tau$. On the other hand, if $g'_\tau=\tau$, then
$ g'_{\tau_n}\le {\tau_n}$ implies that $g'_{\tau_n}\downarrow
g'_{\tau}$.

So whenever $\tau\not=g_\tau$, we have $\rho B_{g_{\tau_n}} = \rho
B_{g_\tau}$ eventually, $L_{\tau_n} \rightarrow L_{\tau}$, and the
latter contains at least two points. By Lemma \ref{fathit} we get
 $$
 \BE(B_{\tau_n}, z, \rho B_{g_\tau},L_{\tau_n}) \cd
 \BE(B_{\tau},   z, \rho B_{g_\tau},L_{\tau}),
 $$
which implies (\ref{cd}).

If $g_\tau=\tau$, then set $K_n:=B[g_{\tau},g_{\tau_n}]$. Note
that $K_n \cup L_{\tau_n} \cup \{B_{\tau_n}\} \subset B[\tau,
\tau_n]$ converges to $\{B_\tau\}$ along a path (i.e. a sub-path
of $B_t$). Therefore Lemma \ref{slimhit} yields (\ref{cd}).
\end{proofof}

\section{Conformal maps, paths, semigroups} \lb{s.conformal}

We begin this section with some notation and standard facts about
complex geometry. Most can be found in \cite{pommerenke},
\cite{lawlersle}, and best of all the upcoming book
\cite{lawlerbook}, which reviews complex geometry with an eye
towards applications in stochastic processes.

If $A\subset \HH$, let  $h_y(A)$ be harmonic measure, or the
hitting distribution of Brownian motion started at $y$ and stopped
at $A\cup \RR$. Let
$$h_\infty(A)=\lim_{y\rightarrow\infty} y\cdot
h_y(A)$$ Define
\begin{eqnarray}\label{capz}
 \capz(A)&=&h_\infty(A),\\
 \label{capo}
 \capo(A)&=&\int_{\overline A} \Im(z) h_\infty(dz).
\end{eqnarray}
If $D$ is a subset of $\HH$ so that $\HH\setminus D$ is bounded,
then the capacity $\capz$ in $D$ can be defined analogously and is
denoted $\capz(D,A)$. Conformal invariance of harmonic measure
implies that if $f$ is a conformal homeomorphism $D\rightarrow
\HH$ satisfying $f'(\infty)=1$, then $\capz(D,A)=\capz(f(A))$. The
quantity $\capo$ will be referred to as half-plane capacity.

If $A$ is bounded, $A \cup \RR$ is connected, and $f:\HH\setminus
A \to \HH$ is a conformal homeomorphism with hydrodynamic
normalization (i.e. $f(z)-z\to 0$ as $z\to\infty$), then we have
 \bel{hyd}
 f(z)=z+\capo(A)/z+O(1/z^2)
 \ee
as $z\rightarrow \infty$. Such maps have the property that for all
$z\in \HH$
 \bel{shiftdown}
  \Im f(z)\le \Im(z).
 \ee
Recall the following standard fact about complex geometry. There
exist a constant $c$ so that if $A\subset \overline \HH$ is
connected, then
 \bel{capdia}
 c^{-1} \diam(A)\le\capz(A)\le c\,\diam(A),
 \ee
if $A\cup \RR$ is connected and $x+iy\in A$, then
 \bel{capim}
 y^2/4 \le \capo(A).
 \ee
Also, if $A,A'$ are disjoint, then we have
 \bel{suba}
 \capo(A\cup A')\le  \capo(A)+\capo(A')
 \ee
This implies that if $A_1,A_2,A$ are disjoint, with $A_1\cup R$
and $\HH\setminus A_1$ connected, then
 \bel{concap}
\capo(A_1 \cup A_2 \cup A)-\capo(A_1 \cup A) \le \capo(A_1\cup
A_2)-\capo(A_1)
 \ee
We apply the conformal map $f:(\HH\setminus A_1,0,\infty)\to
(\HH\setminus A_1,0,\infty)$; then by additivity of $\capo$ the
desired inequality transforms to $\capo(f(A_2\cup A))-\capo(f(A))
\le \capo(f(A_2))$, which follows from \re{suba}.

Let $\Phi$ denote the set of paths $\pi:[0,\tau]\rightarrow
\overline \HH$ which intersect the real line only at time $0$ and
$\pi(\tau)$ is in the boundary of $\rho\varphi$, the infinite
connected component of the complement of the image $\pi[0,\tau]$.

To each $\pi \in \Phi$ we can associate the unique conformal
homeomorphism $f_\pi:\HH \rightarrow \rho_\pi$ with
hydrodynamic normalization. which fixes $\infty$, has
derivative $1$ there, and extends to the boundary $\RR$
continuously (see Remark \ref{curve}) and maps $0$ to $\pi(t)$.
Such a map has expansion
 \bel{shiftexp}
 f(z)=z + a_0 + az^{-1} +
O(z^{-2})
 \ee
at $\infty$. The coefficient $a=\capo(\pi)$ behaves additively
under compositions of maps. It has the scaling property
 \bel{scaling}
 a(r \pi)=r^2 a(\pi).
 \ee
Using the conformal maps, it is possible to compose two paths in
$\Phi$. Let $\pi\circ \pi'$ equal $\pi$ on $[0, \ell_\pi]$, and
equal the $f_\pi$-transform of $\pi'$ (see the definition below)
time-shifted by $\ell_\pi$ for the rest of the time interval.

Thus $\Phi$ is a semigroup, and the map $a$ is a semigroup
homomorphism $\Phi \rightarrow [0,\infty)$.

\begin{definition} \rm Let $\pi$ be a path in a regular domain $D_1$
(possibly starting and ending at $\partial D_1$), and let $f$
be a conformal homeomorphism $D_1 \rightarrow D_2$. If the
integral
 \bel{timechange}
s(t)=\int_0^t |f'\circ \pi(\tau)|^2 d\tau
 \ee
is finite for each $t$, then let $t(s)$ be the inverse of $s(t)$;
the function $f\circ \pi\circ t (s)$ is called the {\bf
$f$-transform} of $\pi(t)$.
\end{definition}

\section{Independent increments and a local time for beads}
\lb{s.independent}

Consider the process $B\sim\BE$, and let $G$ denote the set of
cuttimes. If $B(t)$ is on the boundary of $\rho B_t$, let $f_t$
denote the normalized conformal homeomorphism mapping $(\rho
B_t,B(t), \infty)$ to $(\HH,0,\infty)$. In this section, we
explore the following consequence of Proposition \ref{mp} and the
conformal invariance of Brownian excursion.

\begin{corollary}[Independent increments at cuttimes]
\label{iip} Let $\tau$ be a ${\mathcal G}_t$-stopping time
supported on the set of cuttimes $G$. Then the $f_\tau$-map of the
future has distribution $\BE$.
\end{corollary}

By a minor abuse of notation,  let $a(t)$ denote the half-plane
capacity of $B[0,t]$.

\begin{lemma} The set $a(G)$ is a closed regenerative set
with a scale-invariant distribution.
\end{lemma}
Recall that a random set $S$ has scale-invariant distribution for
every $c>0$ the sets $cS$ are identically distributed. It is
called regenerative if for every $\{\sigma[S \cap [0,t]];\; t\ge
0\}$-stopping time $\nu$ which satisfies $\nu\in S$ a.s. the
translated set $(S-\tau) \cap [0,\infty)$ has the same
distribution as $S$.

\begin{proof}
Closedness follows from the continuity of half-plane capacity and
since the set of cuttimes is closed. Scale-invariance follows from
the scale-invariance of Brownian excursion and the scaling
property $\re{scaling}$ of the half-plane capacity $a$.

To check the regenerative property, assume that the stopping time
$\nu$ is as above. Then $\tau=a^{-1}(\nu)$ is a stopping time for
the cutpoint filtration $\GG$, and $\tau=a^{-1}(\nu)$ is supported
on cutpoints. Therefore by Corollary \ref{iip} the distribution of
the $f_\tau$-mapping of the future of $B$ after $\tau$ has
distribution $\BE$.

In particular, since $a$ behaves additively under conformal
homeomorphisms, the distribution of $(a(G)-\nu)\cap[0,\infty)$
given the past of the excursion is the same as the distribution
of $a(G)$. Thus $a(G)$ is a regenerative set, as required.
\end{proof}

By a result of \cite{kingman73}, the above implies that $a(G)$ is
the image of a stable subordinator. More precisely, there exists a
nondecreasing random function $(\lambda(a),\, a \ge 0)$ adapted to
the filtration $(\sigma(a(G) \cap [0,a]),\, a \ge 0)$ so that
$\lambda$ increases exactly on the set $a(G)$ its right-continuous
inverse $(a_\lambda,\, \lambda\ge 0)$ is a stable subordinator
with index $\alpha\in (0,1)$. The normalization can be chosen so
that the L\'evy measure (the intensity of the Poisson point
process of jumps) assigns mass $y^{-\alpha}$ to the interval
$(y,\infty)$ for all $y\ge 0$. We call $\lambda$ the {\bf bead
(local) time}.

Let $B_{\rightarrow \lambda}$ denote the process $(B(t), 0\le t
\le t(\lambda))$ as an element of the semigroup $\Phi$. Recall
that a L\'evy process $X$ on a topological semigroup $\Phi$ is a
right continuous process with left limits with the property that
at any fixed time $t$, given the entire past $(X_s, 0\le s \le
t)$, the distribution of $X$ is the same as the past composed with
an independent copy of $X$. Since the inverse bead local times are
$\GG_t$-stopping times supported on $G$, Corollary \ref{iip}
immediately gives

\begin{proposition}$(B_{\rightarrow \lambda},\, \lambda\ge 0)$ is a
$\Phi$-valued L\'evy process.
\end{proposition}

Let $t(\lambda)$ be the time corresponding to local time
$\lambda$, i.e. the solution of $a(t(\lambda))=a_\lambda$.
Whenever the process $a_\lambda$ has a jump, the segment
$(B(t),\,t\in[t(\lambda-),t(\lambda+)])$ has no cutpoints, while
$t(\lambda-),\,t(\lambda+)$ are cuttimes. Let $\beta(\lambda)$
denote the $f$-mapping of this segment from $\rho B(t(\lambda-))$
to $\HH$. Then $\beta(\lambda)$ is an element of the semigroup
$\Phi$. When $a(\lambda)$ has no jump, we set $\beta(\lambda)$
equal some null state.

Note the following deterministic fact. For an interval
$I=[\lambda_1,\lambda_2]$ the set $\beta(I)$ equals the set
$\beta'(I-\lambda_1)$ for the process which is the
$f_{a(\lambda_1)}$-image of the original one. Since the inverse
bead local times are stopping times, the independent increment
property of Corollary \ref{iip} implies that for non-overlapping
intervals $I$, the sets $\beta(I)$ are independent. We have shown

\begin{proposition} $\beta$ is a
$\Phi$-valued Poisson point process.
\end{proposition}

The {\bf Brownian bead measure} ($\BB$) is the $\sigma$-finite
intensity measure of the Poisson point process $\beta$. In simple
terms, for a set of paths $A$ the measure $\BB(A)$ equals the
expected number of elements of $A$ among the beads $\beta[0,1]$.

\section{Properties of beads}\lb{s.properties}

The goal of this section is to establish some simple properties of
the {\bf bead process}, i.e. the measure $\BB$.

{\bf Scaling.} It follows from the scaling properties of $\BE$ and
half-plane capacity that for $c>0$ we have
$\BB(rA)=r^{-2a}\BB(A)$, where scaling a path by $c$ means scaling
in space by a factor $c$ and time accordingly. This implies that
the law $\BB$ can be decomposed as a product of measures on
``shape space'' and ``size space''. Let $\beta$ be chosen
according to $\BB$ conditioned to have size $a$ at least $1$, and
rescale $\beta$ to have size 1. The resulting probability measure
$\nbb$, determines the shape of $\BB$, while the size is given
independently by the measure $d(y^{-\alpha})$.

We now check that the Markov property for cuttimes and conformal
invariance imply a Markov property of beads.

Note that the time of beads is only defined up to translation
(since the integral in the $f$-transform may not be finite). We
may pick a rule to set $t=t_1$ for some fixed number $t_1$ when
(and if) the bead grows to size 1 (which here is an arbitrary
positive number). Let $t_0$ be the starting time of the bead
(possibly $-\infty$), and let $A_t$ denote the image of the time
interval between $t_0$ and the first $t$-local cuttime after $t_0$
under $B$. Let $T\ge t_1$ be a stopping time with respect to the
canonical filtration generated by the past of the process.

\begin{lemma}[Markov property of beads]\lb{mpb}
Under the measure $\BB$, the process $(B(T+s),\, s\ge 0)$ has the
same distribution as a process $B'$ distributed as
$\BE(B(T),\infty, \HH,A_T)$ and stopped at a random time $\tau$.
Here $\tau$ is the first time that
$$(B[t_0,\tau]\cup B'[0,\tau))\cap B'(\tau,\infty] =0.$$
\end{lemma}
\begin{proof}
It suffices to show this for the case $T=t_1$ (recall that the
bead starts at a time $t_0<t_1$). The general case follows from
first applying the lemma at $T=t_1$ and then using the Markov
property of Brownian excursion. More precisely, if $B'$ is as
given in the claim, and $T\ge t_1$ is a stopping time, then the
distribution of the future of $B$ after $T$ is just the
distribution of the future of $B'$ after $T$, and $B'\sim
\BE(B(T),\infty, \HH,A_T)$ stopped at $\tau$.

For the $T=t_1$ case note that $\BB$ conditioned to have size at
least $1$ is by definition the following. $\BE$ is run until the
first time $S$ that $a(S)-a(g'(S))=1$ where $g'(S), g(S)$ are the
global cuttimes immediately before and after time $S$. Then $\BB$
is the $f_{g'(S)}$-map of $B(g'(S)+s,\,0\le s \le g(S)-g'(S))$.

Therefore the statement of the lemma follows from the Markov
property of Proposition \ref{mp}, conformal invariance and the
fact that $f$-mappings preserve half-plane capacity.
\end{proof}

Note that an example of such a stopping time is the first time
when $\BB$ hits the line with imaginary part $y$. This will be
used in Lemma \ref{lifetime} to show that in fact beads have
finite lifetime, and can be started at $t_0=0$. Once this has been
done, it is straightforward to extend Lemma \ref{mpb} to arbitrary
stopping times $T>t_0$.

\begin{remark}[Beads determine the excursion]\lb{determine}\rm Another property shared by
beads and It\^o excursions is that the process $(\beta_\lambda,\,
\lambda \ge 0)$ determines the process $(B(t),\, t\ge 0)$. When
$B(t)$ is on the boundary of $B[0,t]$ (these are the so-called
pioneer points), the parameter $a_0$ in the conformal shift
\re{shiftexp} defines the ``horizontal'' location of $B(t)$. Given
$\beta$, it is straightforward to determine $a_0$ as a function of
the half-plane capacity $a$ for each $a\ge 0$. This gives a
L\"owner chain for the pioneer points of $B$ (see, for example,
\cite{lawlersle} for definitions), and hence determines its outer
boundary. Now assume that the outer boundaries of $B$ and $B'$
agree, but they still differ in some way within a particular bead.
Since $f$-mappings are one-to-one, then that bead has to be mapped
by the same conformal shift to different points in the processes
$\beta, \beta'$. This answers a question posed by an anonymous
referee.
\end{remark}

The proof of the following simple fact is left to the reader.

\begin{lemma}  \lb{time}
Let $A$ be a nonempty subset of $\RR\times (0,1]\subset \HH$ so
that $A\cup \RR$ is connected. Let $B$ have distribution
$\BE(z,\infty, \HH)$ conditioned to hit $A$ and stopped when this
happens. Let $T_1$ be the time $B$ spends in the strip with
imaginary parts between $1$ and $2$. There exists absolute
constants $c,\gamma$ so that for all $t>0$ we have
$$
\pr(T>t) < c e^{-\gamma t}.
$$
\end{lemma}

\begin{lemma}[Finite lifetime]\lb{lifetime}
Let $T$ denote the lifetime of the process with distribution
$\BB$. Then a.e. $T<\infty$.
\end{lemma}

\begin{proof}
Let $B$ be a bead conditioned to hit $\Im z=y$. By the Markov
property, the future of $B$ after this time is just $\BE(z,\infty,
\HH)$ conditioned to hit a certain subset of the past. Let
$T(y,2y)$ denote the time $B$ spends with imaginary part in this
interval. By Lemma $\ref{time}$ and scale invariance
$$
\pr(T(y,2y)>ty^2) < ce^{-\gamma t}
$$
and therefore, for the unconditioned measure
$$
\BB(T(y,2y)>ty^2) < y^{-\alpha}c'e^{-\gamma t}.
$$
setting $y_n=2^{-n}$, and $t_n=c_1 n$ the right hand side becomes
$ 2^{\alpha n}2^{-c_1\gamma n} $ which is summable for an
appropriate choice of $c_1$. By the first Borel-Cantelli lemma
(which also holds for $\sigma$-finite measures) we get that $\BB$
a.e. for all large $n$
$$
T(2^{-n},2^{1-n}) \le c n 2^{-2n}
$$
summing this we get that for some ``random'' constant $K$ and all
$y\le 1$ \bel{timeineq} T(0,y)< K y^2|\mbox{log}\,y| \ee But the
Markov property and the existence of cutpoints for large times
implies that after hitting $\Im(z)=1$, $B$ will have finite
lifetime $T'$ a.s. Therefore $T\le T(0,1)+T'$ is finite $\BB$-a.e.
as required.
\end{proof}

\section{The exponent giving the bead index }\lb{s.exponent}

The goal of this section is to identify the index $\alpha$ of the
stable process driving Brownian beads as an exponent for a large
deviation event. Let $A(t)$ denote the event that the half-plane
excursion $B$ has no cuttime between times $1$ and $t$.

\newtheorem*{exponent.thm}{Theorem \ref{exponent}}
 \begin{exponent.thm}
For large $t$, we have $\pr A(t) = t^{-\alpha + o(1)}$.
\end{exponent.thm}

\cite{math.PR/0302115} recently computed essentially the same
exponent using generalized SLE processes, with the result
$\alpha=1/2$. We will compute $\alpha$ directly by a simpler
argument in the next section.

Let $(a(\lambda),\lambda \ge 0)$ be a stable subordinator with
index $\alpha$. We will use the following simple fact. It follows
directly from \cite{bertoin} page 76 Proposition 2.

\begin{fact}\label{bert} Let
\[
X =\min\left(a([0,\infty)) \cap [1,\infty) \right)
\]
then for some positive $c$ as $x\rightarrow \infty$ we have $
\pr(X>x) \sim c x^{-\alpha}. $
\end{fact}

Let $A'(a)$ denote the event that there is no cuttime so that the
past has half-plane capacity between $1$ and $a$. Fact \ref{bert}
implies
\begin{equation}\label{alpha}
\pr A'(a) \sim c a^{-\alpha}\quad\mbox{ as } a \rightarrow \infty.
\end{equation}
In order to conclude Theorem \ref{exponent}, we only need to show
that half-plane capacity and time are not too far from each other;
in fact it suffices to show the following.

\begin{lemma}
We have $\pr(t/s<\capo\,B[0,t]<ts) \ge 1-ce^{-\gamma s}$.
\end{lemma}

\begin{proof}
For a set $A$ in the plane, let $M_x$ and $M_y$ denote the sup of
the absolute value of the projection of $A$ to the $x$ and $y$
axes, respectively. If $A\subset \overline \HH$ contains zero and
is connected, then by (\ref{capim}) and considering the half-plane
capacity of a rectangle we get
$$
M_y^2/4 \le \capo(A) \le c\max(M_y^2,M_x).
$$
Now let $A=B[0,1]$. Then it is easy to check the following simple
property of the maxima of Brownian motion and the 3-dimensional
Bessel processes.
$$
\pr( 1/s<M_y^2  \le \max(M_y^2,M_x)<s) \ge 1-ce^{-\gamma s},
$$
and the claim follows by scale-invariance.
\end{proof}

\section{The value of the bead index $\alpha$}
\lb{s.index}

Consider the process $B\sim \BE$; in this section we will index
$B$ either by time or bead local time $\lambda$; in the latter
case we will stick to the notation $\lambda$. Let $\mathcal L$
denote Lebesgue measure, let $\mu$ denote the random measure on
the half plane given by
$$
\mu(A)={\mathcal L}(\lambda\,:\,B(\lambda)\in A),
$$
and let $\overline \mu(A):=\ev \mu(A)$. By the scaling of
half-plane capacity $a(\lambda)$ and the scale-invariance of
Brownian motion we get that for $r\ge 0$
$$
\omu(rA) = r^{2\alpha} \omu(A).
$$

\begin{lemma} \lb{strip}
The $\omu$-measure of $\,\Strip_1=\RR \times (0,1]$
is finite.
\end{lemma}
\begin{proof}
Consider the random measure
$$
\mu(\Lambda,A)= {\mathcal L}(\lambda\in \Lambda\,:\,B(\lambda)\in
A)
$$
and let $\omu(\Lambda,A):=\ev\mu(\Lambda,A)$. Let $\lambda$ be the
first bead time at least $1$ so that $B(\lambda)\in \Strip_1$, and
let $f$ denote the corresponding conformal shift. Let
$p=\pr(\lambda<\infty)$; it is easy to check that $p<1$. For $n\ge
0$ we have
\begin{eqnarray*}
\ev(\mu([1,n+1],\Strip_1) \gi \GG_\lambda) &\le&
\one(\lambda<\infty) \omu([0,n+1-\lambda],f(\Strip_1))\\&\le&
\one(\lambda<\infty)\omu([0,n],\Strip_1)
\end{eqnarray*}
because of the cuttime Markov property (Proposition \ref {mp}) and
the fact \re{shiftdown} that $f(\Strip_1)\subset \Strip_1$.
Taking expected values gives
$$\omu([1,n+1],\Strip_1) \le p
\omu([0,n],\Strip_1)=:p \omu_n.$$ This yields the recursion
$\omu_{n+1} \le 1 +p \omu_n$, which gives the bound $\omu_\infty
\le 1/(1-p)$, as required.
\end{proof}

\begin{lemma}\lb{absolute}
$\omu$ is absolutely continuous with respect to Lebesgue measure
on $\HH$ with density bounded below and above on compacts.
\end{lemma}

\begin{proof}
Let $A\subset \HH\setminus\{0\}$ be a closed set so that $A\cup
\RR$ and $\HH \setminus A$ are connected, and let $f$ be the
conformal homeomorphism $(\HH\setminus A,0,\infty)\to
(\HH,0,\infty)$. Call such $f$ a subdomain map.

Fix a path $B$ which avoids $A$, and consider its image $f(B)$.
The cuttimes $g$ of $B$ can be parameterized by $a(g)$ i.e.
$\capo$ of the past, as well as $a'(g)$ that is $\capo$ of the
past of the image $f(B)$. Fact \re{concap} implies that if
$g_1<g_2$ are two cuttimes, then
$$
a(g_2)-a(g_1) \le a'(g_2)-a'(g_1),
$$
in particular, beads are smaller when measured by $a'$ then when
measured by $a$. Let $E\subset \HH \setminus A$ a generic Borel
subset, and let $\mu'(E)$ denote $\mu(E)$ measured for the
$f$-mapping of the process $B$. Since a.s. $\mu(E)$ can be
computed as the $\eps\to 0$ limit of the rescaled number of beads
of size at least $\eps$ starting at a cutpoint in $E$, it follows
that we have $\mu'(f(E))\le \mu(E)$. By the restriction property,
$$
\omu(f(E)) = \ev\left(\mu'(f(E))|B\mbox{ avoids }A\right) \le
\ev(\mu(E)|B\mbox{ avoids }A),
$$
and therefore
$$
\omu(f(E)) \le \ev(\mu(E))/\pr(B\mbox{ avoids }A) =
\omu(E)(f'(0)f'(\infty))^{-1}.
$$
Let $D_r(z)$ denote the open disk of radius $r$ about $z$. Let $K$
be a compact subset of $\HH$. If $z\in K$ and $w$ is sufficiently
close to $0$, then there exists a subdomain map $f$ so that
$f(z)=w$. Consider the set of points $w$ for which $f$ exists for
all $z\in K$; this set contains a rescaled version $sK$ of $K$. By
compactness, we may choose subset maps for each $z\in K,sw\in sK$
so that
\begin{eqnarray*}
c_0&\le& f'(0)f'(\infty),  \\
D_{sr}(sw)&\subseteq &f(D_{c_2r}(z))\qquad\mbox{ for all }r<c_1
\end{eqnarray*}
with constants $c_0,c_1,c_2$ depending on $K$ only. Then
\begin{eqnarray*}
s^{2\alpha}\omu(D_{r}(w)) &=& \omu(D_{sr}(sw)) \\&\le&
\omu(f(D_{c_2r}(z))) \\&\le&
(f'(0)f'(\infty))^{-1}\omu(D_{c_2r}(z))\\&\le&
c_0^{-1}\omu(D_{c_2r}(z)).
\end{eqnarray*}
This uniform bound implies the claim by standard arguments about
absolute continuity.
\end{proof}

The following are needed for the proof of
 \newtheorem*{indexhalf.thm}{Theorem \ref{indexhalf}}
 \begin{indexhalf.thm}
 \indexhalfthm
 \end{indexhalf.thm}

Let $D_r=D_r(i)$ denote the open disk of radius $r$ about $i$. Let
$\tau$ be the (possibly infinite) first hitting time of $D_r$ for
the process $B\sim\BE$. Let $B^{(1)}$ denote the path image
$B[0,\tau]$. Recall the definition of the capacity $\capz(D,A)$ of
$A$ from $\infty$ in $D$ from Section \ref{s.conformal}, and use
the shorthand $\capz(t,A)= \capz(\rho B_t,A)$. Let
$K_r=\capz(\tau,D_r)$, and recall the notation $x_r \asymp y_r$
for the existence of a constant $c>0$ so that $c^{-1}y_r \le
x_r\le c y_r$ (here for all small $r$).

\begin{proposition}\label{zbounds}
As $r\to 0$ we have $ \omu(D_r) \asymp |\mbox{\rm log } r|\ev
K_r^{1+2\alpha}. $
\end{proposition}

Let $B^{(2)}$ denote the path after the last exit from $D_r$. Let
$Z'_r=\pr (B^{(1)} \nis B^{(2)}\,|\,\FF_\tau)$, where ``$\nis$''
denotes ``does not intersect''. Let $Z_r$ denote the probability
given $\FF_\tau$ that the image of an independent process
$\BE(\infty, D_r,\HH)$ does not intersect $B^{(1)}$; recall that
this image can be defined via conformal mapping of $\BE$ started
at a more conventional boundary point.

\begin{lemma}\lb{sizes} We have
$|\mbox{\rm log}\ r| K_r \asymp Z_r \asymp Z'_r$, with
deterministic constants.
\end{lemma}

\begin{proof}
We will use the time-reversal property of Brownian excursion: if
$B\sim \BE(a,z,D)$ and ends at random time $\tau$, then
$(B(\tau-t),\, t\in[0,\tau])\sim \BE(z,a,D)$. This, as well as
conformal invariance, can be used to define the image of
$\BE(\infty,z,\HH)$ (together with a time-parameterization
starting at $-\infty$, but we won't need this).

The quantities in question are given by the measures of paths that
do not intersect $B^{(1)}$ before hitting $D_r$ under the measures
$|\mbox{log}\ r|\BE(\infty,\RR,\HH)$, $\BE(\infty, D_r,\HH)$, and
$\BE(\infty,B(\tau),\HH)$, respectively. It follows from the
definition of $\BE$ that the distribution of these paths up to the
hitting time of $D_r$ given the hitting position agree. It is easy
to check that the in each case hitting position has a smooth
density bounded below and above with respect to uniform measure on
$\partial D_r$. The $|\mbox{log}\ r|$ normalizing factor is
necessary so that the measure of paths hitting $D_r$ is bounded
below and above by constants as $r\to 0$.
\end{proof}

\begin{proofof}{Theorem \ref{indexhalf}}
By the Lemma \ref{absolute}, as $r\rightarrow 0$ we have
$$
\omu(D_r) \asymp r^2.
$$
We have
$$
\ev K_r^\gamma \asymp |\mbox{log }r|^{-\gamma} \ev Z_r^\gamma =
r^{\xi(1,\gamma)+o(1)},$$ where the first approximation follows
from Lemma \ref{sizes}, and the second is one definition of the
intersection exponent. More precisely, we should consider the
analogue of $Z_r$ for $B^{(1)'}$, which has distribution
$\BE(0,D_r,\HH)$, but this is absolutely continuous with density
bounded above and below with respect to the distribution of
$B^{(1)}$ given that it hits $D_r$ (see the proof of Lemma
\ref{sizes}).

Thus by Proposition \ref{zbounds} we get $
r^{\xi(1,1+2\alpha)+o(1)} = r^2,$ and the theorem follows by the
monotonicity of $\xi(1,\cdot)$ and the known value $\xi(1,2)=2$
(see \cite{lawler95} for these properties of $\xi$).
\end{proofof}

\begin{proofof}{Proposition \ref{zbounds}}\ \\
\indent{\bf Upper bound.} We write
$$
\ev(\mu(D_r)\,|\,\FF_\tau) = \ev(\mu(D_r) \,|\, B^{(1)}\nis
B^{(2)}, \FF_\tau) \pr (B^{(1)} \nis B^{(2)}\,|\,\FF_\tau).
$$
Let $f$ denote the conformal shift at the first cuttime $g$ after
time $\tau$. We have
$$K_r=\capz(\tau,D_r)\ge \capz(g,D_r)=\capz(f(D_r)).
$$
Here, by slight abuse of notation, $f(D_r)$ denotes the image of
the part of $D_r$ on which $f$ is defined. By \re{capdia} the
above implies that $f(D_r)\subset \RR\times (0,cK_r].$ By scaling
and Lemma \ref{strip} we have $\omu(\RR\times
(0,K_r/2])=cK_r^{2\alpha}.$ Thus by the cuttime Markov property
(Proposition \ref{mp}) we get
$$
\ev(\mu(D_r) \,|\, B^{(1)}\nis B^{(2)}, \FF_\tau) =
\ev(\omu(f(D_r))\,|\,\FF_\tau) \le cK_r^{2\alpha}.
$$
Here and in the sequel $c$ denotes a constant whose value may
change by line to line. By Lemma \ref{sizes}
$$
\ev(\mu(D_r) \gi \FF_\tau) \le c Z'_rK_r^{2\alpha} \le c|\mbox{log
}r| K_r^{1+2\alpha},
$$
which implies the upper bound in Proposition \ref{zbounds}.

\bigskip
{\bf Lower bound.} Let $A_1$ denote the event that $B[0,\tau]$
does not intersect the set $ \{w\in D_{2r}\,:\,\arg(w-i)\in
(\pi/4,7\pi/4)\}.$ It is easy to check (see, for example
\cite{lwssharp}) that for all $\gamma\in [1,3]$ and $r<1/2$ we
have
 \bel{za1}
 \ev(Z_r^\gamma;\,A_1) \ge c\ev Z_r^\gamma.
 \ee
Consider the event
$$
A=\{B(g)\in D_{r},\, B[0,g]\nis D_{r/2},\, \capz(g,D_{r/32})\ge
cK_r \}.
$$
This implies $B^{(1)}\nis B^{(2)}$, and it is easy to check that
there is an absolute constant  $c_1$ so that
$$
\pr(A\,|\,B^{(1)}\nis B^{(2)},\FF_\tau) \ge c_1 \one_{A_1}.
$$
We have
\begin{eqnarray*}
 \ev( \mu(D_{r}) \gi \FF_\tau) &\ge&
 \ev(\mu(D_{r});\, A \gi \FF_\tau)\one_{A_1} \\&\ge&
 \ev(\mu(D_{r}) \gi A,\FF_\tau)
 \pr(A \gi B^{(1)}\nis B^{(2)}, \FF_\tau)
 \pr(B^{(1)}\nis B^{(2)}\,|\, \FF_\tau)\one_{A_1}.
\end{eqnarray*}

By the cuttime Markov property (Proposition \ref{mp}) we get
$$
 \ev(\mu(D_{r}) \gi A,\FF_\tau)
= \ev(\omu(f(D_{r})) \gi A,\FF_\tau).
$$
The event $A$ implies that the domain of $f$ contains $D_{r/2}$,
and therefore by the K\"obe quarter theorem $f(D_{r/2})$ contains
a ball $D'$ of radius $|f'(i)|r/8$ centered at $f(i)$. By the same
theorem applied to $f^{-1}$, we get that $f^{-1}(D')$ contains
$D_{r/32}$. By monotonicity and conformal invariance of harmonic
measure, we get
$$
\mbox{diam}(D') \ge c\capz(D')=c\capz(g,f^{-1}(D')) \ge
c\capz(g,D_{r/32}) \ge c' K_r,
$$
where the last inequality is an assumption in $A$. Since
$f(\partial D_r)$ contains an curve $K$ that separates $0$ and
$D'$ from $\infty$ in $\HH$, we have
$$
\dist(D',0)\le \mbox{diam}(K)\le c\,\capz(K)\le
c\,\capz(f(D_r))\le cK_r.
$$
When $K_r=1$, this implies that $\omu(D')$ is bounded below by a
constant. The scaling property then implies that in general
$\omu(D')\ge cK_r^{2\alpha}$. Putting all the above together and
using Lemma \ref{sizes} we get
$$ \ev (\mu(D_{r})\,|\,\FF_\tau) \ge c
Z'_r K_r^{2\alpha}\one_{A_1} \ge c  |\mbox{log }
r|K_r^{1+2\alpha}\one_{A_1}
$$ and the bound follows from (\ref{za1}).
\end{proofof}
\subsection*{Open questions and conjectures}

\begin{question}\rm
The parameter $a_0$ in the conformal shift \re{shiftexp} defines
the ``horizontal'' location of a cutpoint. It follows from our
proofs that the process $((a(\lambda),a_0(\lambda)),\, \lambda \ge
0)$ is a 2-dimensional L\'evy process stable under scaling with
exponents $1/2$ and $1$ in the two respective coordinates. In
particular $(a_0(\lambda),\, \lambda>0)$ is a Cauchy process. What
is the joint distribution of the two processes? (This question may
be better considered in the $\mbox{SLE}_6$ framework.)
\end{question}

\begin{question}\rm
Is it possible to interpret bead local time as a local time at
zero of some process? This may be better answered by considering a
version of $\mbox{SLE}_6$.
\end{question}

\begin{conjecture}\rm
Consider the $\sigma$-finite bead process conditioned to survive
up to distance $r$ (or time $r^2$, or imaginary part $r$), and let
$r\rightarrow \infty$. The limiting process exists and has the
restriction property. The exponent (Werner, personal
communication) should equal 2 (see \cite{lwsrest} for
definitions).
\end{conjecture}

\begin{question}\rm
The bead time $\lambda$ defines a random measure on the
half-plane. The closed support of this measure is the set of
cutpoints. Is it possible to derive the Hausdorff dimension of the
set of cutpoints using this measure, or more generally, using
beads? This would give a new, conceptual proof for the value of
the intersection exponent $\xi(1,1)=5/4$.
\end{question}

\bigskip
\noindent {\bf Acknowledgments.} The author thanks Jim Pitman for
stimulating remarks and recommending the reference \cite{gp80},
which can also be used to construct bead local time. He also
thanks Oded Schramm and Wendelin Werner for inspiring discussions
and conscientious referee for several important remarks and
corrections.
\bibliography{bbb}

\begin{thebibliography}{18}
\expandafter\ifx\csname natexlab\endcsname\relax\def\natexlab#1{#1}\fi
\expandafter\ifx\csname url\endcsname\relax
  \def\url#1{{\tt #1}}\fi

\bibitem[Beffara(2003)]{beffara}
V.~Beffara (2003).
\newblock {Hausdorff Dimensions for $SLE_6$}.
\newblock { arXiv:math.PR/0204208}.
\newblock Preprint.

\bibitem[Bertoin(1996)]{bertoin}
J.~Bertoin.
\newblock {\em L\'evy processes}, volume 121 of {\em Cambridge Tracts in
  Mathematics}.
\newblock Cambridge University Press, Cambridge, 1996.

\bibitem[Burdzy(1989)]{burdzy89}
K.~Burdzy (1989).
\newblock Cut points on {B}rownian paths.
\newblock {\em Ann. Probab.}, 17\penalty0 (3):\penalty0 1012--1036.

\bibitem[Chung(1984)]{MR86d:60088}
K.~L. Chung (1984).
\newblock The lifetime of conditional {B}rownian motion in the plane.
\newblock {\em Ann. Inst. H. Poincar\'e Probab. Statist.}, 20\penalty0
  (4):\penalty0 349--351.

\bibitem[Cranston and McConnell(1983)]{cm83}
M.~Cranston and T.~R. McConnell (1983).
\newblock The lifetime of conditioned {B}rownian motion.
\newblock {\em Z. Wahrsch. Verw. Gebiete}, 65\penalty0 (1):\penalty0 1--11.

\bibitem[Dubedat(2003)]{dubedat}
J.~Dubedat (2003).
\newblock {$SLE(\kappa,\rho)$ martingales and duality}.
\newblock { arXiv:math.PR/0303128}.
\newblock Preprint.

\bibitem[Greenwood and Pitman(1980)]{gp80}
P.~Greenwood and J.~Pitman (1980).
\newblock Construction of local time and {P}oisson point processes from nested
  arrays.
\newblock {\em J. London Math. Soc. (2)}, 22\penalty0 (1):\penalty0 182--192.

\bibitem[Kingman(1973)]{kingman73}
J.~F.~C. Kingman (1973).
\newblock Homecomings of {M}arkov processes.
\newblock {\em Advances in Appl. Probability}, 5:\penalty0 66--102.

\bibitem[Lawler(1995)]{lawler95}
G.~F. Lawler (1995).
\newblock Nonintersecting planar {B}rownian motions.
\newblock {\em Math. Phys. Electron. J.}, 1:\penalty0 Paper 4, approx.\ 35 pp.

\bibitem[Lawler(1996)]{lawler96}
G.~F. Lawler (1996).
\newblock Hausdorff dimension of cut points for {B}rownian motion.
\newblock {\em Electron. J. Probab.}, 1:\penalty0 no.\ 2, approx.\ 20 pp.

\bibitem[Lawler(2001)]{lawlersle}
G.~F. Lawler (2001).
\newblock {An introduction to the stochastic Loewner evolution}.
\newblock Preprint.

\bibitem[Lawler(2003)]{lawlerbook}
G.~F. Lawler (2003).
\newblock Book in preparation.

\bibitem[Lawler et~al.(2002)Lawler, Schramm, and Werner]{lwssharp}
G.~F. Lawler, O.~Schramm, and W.~Werner.
\newblock Sharp estimates for {B}rownian non-intersection probabilities.
\newblock In {\em In and out of equilibrium (Mambucaba, 2000)}, volume~51 of
  {\em Progr. Probab.}, pages 113--131. Birkh\"auser Boston, Boston, MA, 2002.

\bibitem[Lawler et~al.(2003)Lawler, Schramm, and Werner]{lwsrest}
G.~F. Lawler, O.~Schramm, and W.~Werner (2003).
\newblock {Conformal restriction: the chordal case}.
\newblock { arXiv:math.PR/0209343}.
\newblock Preprint.

\bibitem[Lawler and Werner(2003)]{loopsoup}
G.~F. Lawler and W.~Werner (2003).
\newblock {The Brownian loop soup}.
\newblock { arXiv:math.PR/0304419}.
\newblock Preprint.

\bibitem[L{\'e}vy(1940)]{levy40}
P.~L{\'e}vy (1940).
\newblock Le mouvement brownien plan.
\newblock {\em Amer. J. Math.}, 62:\penalty0 487--550.

\bibitem[Pommerenke(1975)]{pommerenke}
C.~Pommerenke.
\newblock {\em Univalent functions}.
\newblock Vandenhoeck \& Ruprecht, G\"ottingen, 1975.
\newblock With a chapter on quadratic differentials by Gerd Jensen, Studia
  Mathematica/Mathematische Lehrb\"ucher, Band XXV.

\bibitem[Werner(2003)]{math.PR/0302115}
W.~Werner (2003).
\newblock {Girsanov's transformation for SLE($\kappa$,$\rho$) processes,
  intersection exponents and hiding exponents}.
\newblock { arXiv:math.PR/0302115}.
\newblock Preprint.

\end{thebibliography}

\end{document}